\documentclass{llncs}

\begin{document}

\newcommand{\C}{{\cal C}}

\def\thefootnote{\fnsymbol{footnote}}

\newcommand{\strike}{{\tt strike}}
\newcommand{\true}{{\tt true}}
\newcommand{\false}{{\tt false}}

\newcommand{\polylog}{{\rm polylog}}

\newcommand{\dmax}{\Delta}

\def\reals{\mathbb{R}}
\def\naturals{\mathbb{N}}

\newcommand{\eqref}[1]{(\ref{#1})}

\newcommand{\pert}[1]{{\widehat #1}}
\newcommand{\proofend}{\smallskip \hfill\mbox{$\Box$}\\}

\newcommand{\ra}{\rightarrow}
\newcommand{\dd}{\delta}
\newcommand{\eps}{\epsilon}
\newcommand{\beq}{\begin{equation}}
\newcommand{\eeq}{\end{equation}}
\newcommand{\bea}{\begin{eqnarray*}}
\newcommand{\eea}{\end{eqnarray*}}
\newcommand{\sub}{\scriptsize}
\newcommand{\mat}{\left(\!\!\begin{array}{cc}}
\newcommand{\rix}{\end{array}\!\!\right)}
\newcommand{\ord}{{O}}
\newcommand{\bra}{\langle}
\newcommand{\ket}{\rangle}
\newcommand{\p}{\partial}
\newcommand{\Z}{{\Bbb Z}}
\newcommand{\dm}{{\rm d}m}
\newcommand{\ds}{{\rm d}s}
\newcommand{\dt}{{\rm d}t}
\newcommand{\dv}{{\rm d}v}
\newcommand{\dx}{{\rm d}x}
\newcommand{\dy}{{\rm d}y}
\newcommand{\dw}{{\rm d}w}
\newcommand{\db}{{\rm d}b}
\newcommand{\dc}{{\rm d}c}
\newcommand{\du}{{\rm d}u}
\newcommand{\dbeta}{{\rm d}\beta}
\newcommand{\dgamma}{{\rm d}\gamma}
\newcommand{\ddelta}{{\rm d}\delta}
\newcommand{\e}{{\rm e}}
\def\limninf{\lim_{n \rightarrow \infty}}

\newcommand{\dg}{{\rm d}g}
\newcommand{\dz}{{\rm d}z}
\newcommand{\dalpha}{{\rm d}\alpha}
\newcommand{\dda}{\frac{{\rm d}}{\dalpha}}
\renewcommand{\d}{{\delta}}
\newcommand{\amax}{\alpha_{\max}}
\newcommand{\bmax}{\beta_{\max}}
\newcommand{\gamax}{\gamma_{\max}}
\newcommand{\gmax}{g_{\max}}
\newcommand{\za}{\vec{\zeta}}
\newcommand{\zamax}{\vec{\zeta}_{\max}}
\newcommand{\zomax}{\vec{\zeta}^{*}}

\newcommand{\gnae}{\ell_r}

\def\as{a.s.\ }
\def\whp{w.h.p.}
\def\Whp{W.h.p.}

\def\eg{e.g.\ }
\def\ie{i.e.\ }
\newcommand{\fo}[2]{\ensuremath{{F}_{#1}(n,#2)}}
\def\var{{\rm var}}
\def\ex{{\bf E}}

\title{On the 2-colorability of random hypergraphs}

\author{
Dimitris Achlioptas and Cristopher Moore\thanks{Supported
by NSF grant PHY-0071139, the Sandia University Research Program, and
Los Alamos National Laboratory.}}
\institute{
Microsoft Research, Redmond, Washington {\tt optas@microsoft.com} \and
Computer Science Department, University of New Mexico, Albuquerque
and the Santa Fe Institute, Santa Fe, New Mexico
{\tt moore@cs.unm.edu} }

\maketitle

\begin{abstract}
A 2-coloring of a hypergraph is a mapping from its vertices to a
set of two colors such that no edge is monochromatic. Let
$H_k(n,m)$ be a random $k$-uniform hypergraph on $n$ vertices
formed by picking $m$ edges uniformly, independently and with
replacement. It is easy to show that if $r \geq r_c = 2^{k-1} \ln
2 - (\ln 2) /2$, then with high probability $H_k(n,m=rn)$ is not
2-colorable. We complement this observation by proving that if $r
\leq r_c - 1$ then with high probability $H_k(n,m=rn)$ is
2-colorable.
\end{abstract}

\section{Introduction}

For an integer $k\geq 2$, a {\em $k$-uniform hypergraph} $H$ is an
ordered pair $H=(V,E)$, where $V$ is a finite non-empty set,
called the set of {\em vertices\/} of $H$, and $E$ is a family of
distinct $k$-subsets of $V$, called the {\em edges\/} of $H$. For
general hypergraph terminology and background
see~\cite{berge_book}. A {\em $2$-coloring} of a hypergraph
$H=(V,E)$ is a partition of its vertex set $V$ into two (color)
classes so that no edge in $E$ is monochromatic. A hypergraph is
{\em 2-colorable} if it admits a 2-coloring.

The property of 2-colorability was introduced and studied by
Bernstein~\cite{bernstein} in the early 1900s for infinite
hypergraphs. The 2-colorability of finite hypergraphs, also known
as ``Property B'' (a term coined by Erd\H os in reference to
Bernstein), has been studied for about eighty years (\eg
\cite{beck_3_chrom_hyper,Erd63_normat,LLL,KLRH,KS,RadSri00,SSU}).
For $k=2$, \ie for graphs, the problem is well understood since a
graph is 2-colorable if and only if it has no odd cycle. For
$k\geq 3$, though, much less is known and deciding the
\mbox{2-colorability} of $k$-uniform hypergraphs is
NP-complete~\cite{Lov_hyper_complex}.

In this paper we discuss the 2-colorability of random $k$-uniform
hypergraphs for $k \geq 3$. (For the evolution of odd cycles in
random graphs see \cite{FlaKnuPit89}.) Let $H_k(n,m)$ be a random
$k$-uniform hypergraph on $n$ vertices, where the edge set is
formed by selecting uniformly, independently and with replacement
$m$ out of all possible ${n \choose k}$ edges. We will study
asymptotic properties of $H_k(n,m)$ when $k \geq 3$ is arbitrary
but fixed while $n$ tends to infinity. We will say that a
hypergraph property $A$ holds {\em with high probability} (w.h.p.)
in $H_k(n,m)$ if $\lim_{n\ra\infty} \Pr[H_k(n,m)\mbox{ has }
A]=1$. The main question in this setting is:

\begin{center}
As $m$ is increased, when does $H_k(n,m)$ stop being 2-colorable?
\end{center}

It is popular to conjecture that the transition from
2-colorability to non-2-colorability is {\em sharp}. That is, it
is believed that for each $k \geq 3$, there exists a constant
$r_k$ such that if $r<r_k$ then $H_k(n,m=rn)$ is \whp\
2-colorable, but if $r>r_k$ then \whp\  $H_k(n,m=rn)$   is not
2-colorable. Determining $r_k$ is a challenging open problem,
closely related to the satisfiability threshold conjecture for
random $k$-SAT. Although $r_k$ has not been proven to exist, we
will take the liberty of writing $r_k \geq r^*$ to denote that for
$r < r^*$, $H_k(n,rn)$ is 2-colorable \whp\ (and analogously for
$r_k \leq r^*$).

A relatively recent result of Friedgut~\cite{frie} supports this
conjecture as it gives a {\em non-uniform\/} sharp threshold for
hypergraph 2-colorability. Namely, for each $k \geq 3$ there exists a
{\em sequence} $r_k(n)$ such that if $r < r_k(n) - \eps$ then \whp\
$H_k(n,rn)$ is 2-colorable, but if $r > r_k(n) + \eps$ then \whp\
$H_k(n,rn)$ is not 2-colorable. We will find useful the following
immediate corollary of this sharp threshold.
\begin{corollary}\label{frie}
If
$${\liminf_{n \ra \infty}} \Pr[H_k(n,r^* n) \mbox{ is
2-colorable}]>0 \enspace ,$$ then for $r < r^*$, $H_k(n,rn)$ is
2-colorable \whp
\end{corollary}

Alon and Spencer~\cite{AloSpeun} were the first to give bounds on
the potential value of $r_k$. In particular, they  observed that
the expected number of 2-colorings of $H_k(n,m=rn)$ is $o(1)$ if
$2(1-2^{1-k})^r<1$, implying \begin{equation}\label{eq:upperbound}
r_k < 2^{k-1} \ln 2 - \frac{\ln 2}{2} \enspace .
\end{equation}
Their main contribution, though, was providing a lower bound on
$r_k$. Specifically, by applying the Lov\'asz Local Lemma, they
were able to show that if $r=c\,2^k/k^2$ then \whp\ $H_k(n,rn)$ is
2-colorable, for some small constant $c>0$.

In~\cite{hyp2col}, Achlioptas, Kim, Krivelevich and Tetali reduced
the asymptotic gap between the upper and lower bounds
of~\cite{AloSpeun} from order $k^2$ to order $k$. In particular,
they proved that there exists a constant $c>0$ such that if $r
\leq c\,2^k/k$ then a simple, linear-time algorithm \whp\ finds a
2-coloring of $H_k(n,rn)$. Their algorithm was motivated by
algorithms for random $k$-SAT due to Chao and Franco~\cite{ChFrUC}
and Chv\'{a}tal and Reed~\cite{mick}. In fact, those algorithms
give a similar $\mathrm{\Omega}(2^k / k)$ lower bound on the
random $k$-SAT threshold which, like $r_k$, can also be easily
bounded as $\mathrm{O}(2^k)$.

Very recently, the authors eliminated the gap for the random
$k$-SAT threshold, determining its value within a factor of
two~\cite{ksat}. The proof amounts to applying the ``second
moment'' method to the set of satisfying truth assignments whose
complement is also satisfying.  Alternatively, one can think of
this as applying the second moment method to the number of truth
assignments under which every $k$-clause contains at least one
satisfied literal {\em and\/} at least one unsatisfied literal,
\ie which satisfy the formula when interpreted as a random
instance of Not-All-Equal $k$-SAT (NAE $k$-SAT).

Here we extend the techniques of~\cite{ksat} and apply them to
hypergraph 2-colorability. This allows us to determine $r_k$
within a small additive constant.
\begin{theorem}\label{thm:hyp}
For every $\eps > 0$ and all $k \geq k_0(\eps)$,
\[
r_k \geq 2^{k-1} \ln 2 - \frac{\ln 2}{2} - \frac{1+\eps}{2}
\enspace .
\]
\end{theorem}

Our method actually yields an explicit lower bound for $r_k$ for
each value of $k$ as the solution to a simple equation (yet one
without a pretty closed form, hence Theorem~\ref{thm:hyp}). Below
we compare this lower bound to the upper bound
of~\eqref{eq:upperbound} for small values of $k$. The gap
converges to 1/2 rather rapidly.
\begin{table}[h]\label{tab:val}
$$ \begin{array}{c|ccccccccccc}
    k                   &       3       &       4       &       5       &            7       &           9       &           11      &   12\\
\hline
\mbox{Lower bound}      &       \;\;3/2\;     &       \;49/12\;   &       \;9.973\;   &            \;43.432\;  &           \;176.570\; &           \;708.925\; &   \;1418.712\\
\mbox{Upper bound}      &       \;\;2.409   &       5.191   &       10.740  &            44.014  &           177.099 &           709.436 &   1419.219
%\mbox{$\displaystyle{\frac{\ln 2}{2}(2^k-1)}$}   & 2.426 & 5.199 & 10.744 & 21.834 & 44.015  &       88.376  & 177.099 & 354.545 & 709.436 & 1419.219
\end{array}
$$
\caption{Upper and lower bounds for $r_k$}
\end{table}
\section{Second moment and NAE $k$-SAT}

We prove Theorem~\ref{thm:hyp} by applying the following version
of the second moment method (see Exercise 3.6 in~\cite{ramr}) to
the number of 2-colorings of $H_k(n,m=rn)$.
\begin{lemma}\label{lemma:sec}
For any non-negative integer-valued random variable $X$,
\begin{equation}\label{eq:second}
 \Pr[X > 0] \,\ge\, \frac{\ex[X]^2}{\ex[X^2]} \enspace .
\end{equation}
\end{lemma}

In particular, if $X$ is the number of 2-colorings of
$H_k(n,m=rn)$, we will prove that for all $\eps
> 0$ and all $k \geq k_0(\eps)$, if $r = 2^{k-1} \ln 2 - \ln 2/2 -
(1+\eps)/2$ then there exists some constant $C=C(k)$ such that
\[
\ex[X^2] < C \times \ex[X]^2 \enspace .
\]
By Lemma~\ref{lemma:sec}, this implies $\Pr[X>0] = \Pr[H_k(n,rn)
\mbox{ is 2-colorable}] > 1/C$. Theorem~\ref{thm:hyp} follows by
invoking Corollary~\ref{frie}.

This approach parallels the one taken recently by the authors for
random NAE $k$-SAT~\cite{ksat}. Naturally, what differs is the
second-moment calculation which here is prima facie significantly
more involved.

We start our exposition by outlining the NAE $k$-SAT calculation
of~\cite{ksat}. This serves as a warm up for our calculations and
allows us to state a couple of useful lemmata from~\cite{ksat}. We
then proceed to outline the proof of our main result, showing the
parallels with NAE $k$-SAT and reducing the proof of
Theorem~\ref{thm:hyp} to the proof of three independent lemmata.

The first such lemma is specific to hypergraph 2-colorability and
expresses $\ex[X^2]$ as a multinomial sum. The second one is a
general lemma about bounding multinomial sums by a function of
their largest term and is perhaps of independent interest. It
generalizes Lemma~\ref{lem:basic} of~\cite{ksat}, which we state
below. After applying these two lemmata, we are left to maximize a
three-variable function parameterized by $k$ and $r$. This is
analogous to NAE $k$-SAT, except that there we only have to deal
with a one-variable function, similarly parameterized. That
simpler maximization, in fact, amounted to the bulk of the
technical work in~\cite{ksat}. Luckily, here we will be able to
get away with much less work: a convexity argument will allow us
to reduce our three-dimensional optimization precisely to the
optimization in~\cite{ksat}.

\subsection{Proof outline for NAE $k$-SAT}
Let $Y$ be the number of satisfying assignments of a random NAE
$k$-SAT formula with $n$ variables and $m = rn$ clauses. It is
easy to see that $\ex[Y] = 2^n (1-2^{1-k})^{rn}$. Then $\ex[Y^2]$
is the sum, over all ordered pairs of truth assignments, of the
probability that both assignments in the pair are satisfying. It
is not hard to show that if two assignments assign the same value
to $z=\alpha n$ variables, then the probability that both are
satisfying is
\[
p(\alpha) = 1 \,-\, 2^{1-k} \left(2-\alpha^k - (1-\alpha)^k
\right) \enspace .
\]
Since there are $2^n {n \choose z}$ pairs of assignments sharing $z$
variables, we have
\[  \frac{\ex[Y^2]}{\ex[Y]^2}
 =  \sum_{z=0}^n {n \choose z} \left[\frac{1}{2} \,
         \left(\frac{p(z/n)}{(1-2^{1-k})^2}\right)^r\right]^n
\enspace .
\]
To bound such sums within a constant factor, we proved the
following in~\cite{ksat}.
\begin{lemma}\label{lem:basic}
Let $f$ be any real positive analytic  function and let
$$
S = \sum_{z=0}^n {n \choose z} \,f(z/n)^n \enspace .
$$
Define $0^0 \equiv 1$ and let $g$ on $[0,1]$ be
\[
g(\alpha) = \frac{f(\alpha)}
                 {\alpha^\alpha \,(1-\alpha)^{1-\alpha}} \enspace
.
\]
If there exists $\amax \in (0,1)$ such that $g(\amax) \equiv \gmax >
g(\alpha)$ for all $\alpha \neq \amax$, and $g''(\amax) < 0$, then
there exist constants $B$ and $C$ such that for all sufficiently large
$n$
\[
B \times \gmax^n \,\le\, S \le\, C \times \gmax^n \enspace .
\]
\end{lemma}
Thus, using Lemma~\ref{lem:basic}, bounding $\ex[X^2]/\ex[X]^2$
reduces to maximizing
\begin{equation}\label{eq:common}
\gnae(\alpha)
 = \frac{1}{2 \,\alpha^\alpha \,(1-\alpha)^{1-\alpha}} \,
   \left(\frac{p(\alpha)}{(1-2^{1-k})^2}\right)^r \enspace .
\end{equation}

Note now that $\gnae(1/2) = 1$ for all $r$ and that our goal is to
find $r$ such that $\gnae(\alpha) \leq 1$ for all $\alpha \in
[0,1]$. Indeed, in~\cite{ksat} we showed that
\begin{lemma}\label{gleas}{\rm \cite{ksat}}
For every $\eps > 0$, and all $k \geq k_0(\epsilon)$, if
\[
r \leq 2^{k-1} \ln 2 - \frac{\ln 2}{2} - \frac{1+\eps}{2}
\]
then $\gnae(1/2) = 1 > \gnae(\alpha)$ for all $\alpha \neq 1/2$ and
$\gnae''(1/2) < 0$.
\end{lemma}
Thus, for all $r,k,\eps$ as in Lemma~\ref{gleas}, we see that
Lemma~\ref{lem:basic} implies $\ex[Y^2]/\ex[Y]^2 < C \times
\gnae(1/2)^n = C$, concluding the proof.

\section{Proof outline for hypergraph 2-colorability}

Let $X$ be the number of 2-colorings of $H_k(n,rn)$.  Let $ q=
1-2^{1-k} $ and
$$
p(\alpha, \beta, \gamma)
  =  1 - \alpha^k - (1-\alpha)^k - \beta^k - (1-\beta)^k
  + \gamma^k + (\alpha-\gamma)^k
  + (\beta-\gamma)^k + (1-\alpha-\beta+\gamma)^k.
$$

We will prove that
\begin{lemma}\label{lem:prob_calc}
There exists a constant $A$ such that
\[ \frac{\ex[X^2]}{\ex[X]^2} \leq
   \frac{1}{A^2} \sum_{z_1+\cdots+z_4=n}
   {n \choose z_1,z_2,z_3,z_4} \,
   \left( \frac{1}{4}
   \left( \frac{
   p\!\left(\frac{z_1+z_2}{n},
            \frac{z_1+z_3}{n},
            \frac{z_1}{n}\right) }{q^2} \right)^r \right)^n
\enspace .
\]
\end{lemma}

Similarly to NAE $k$-SAT we would like to bound this sum by a
function of its maximum term.  To do this we will establish a
multidimensional generalization of the upper bound of
Lemma~\ref{lem:basic}.

\begin{lemma}\label{lem:cris}
Let $f$ be any real positive analytic function and let
\[
S \;\;= \sum_{z_1+\cdots+z_{d} = n}{{n \choose z_1, \cdots, z_d}
        \,f(z_1/n,\cdots,z_{d-1}/n)^n}
\enspace .
\]
Let $\displaystyle{ Z = \left\{(\zeta_1,\ldots,\zeta_{d-1}) :
\zeta_i \geq 0 \mbox{ for all $i$, and }\sum%_{i=1}^{d-1}
\zeta_i \le 1\right\}}$. Define $g$ on $Z$ as
\[
g(\zeta_1,\ldots,\zeta_{d-1})
   = \frac{f(\zeta_1,\ldots,\zeta_{d-1})}
          {\zeta_1^{\zeta_1} \cdots
       \zeta_{d-1}^{\zeta_{d-1}} \,
       (1-\zeta_1-\cdots-\zeta_{d-1})^{1-\zeta_1-\cdots-\zeta_{d-1}}}
\enspace .
\]
If i) there exists $\zamax$ in the interior of $Z$ such that for
all $\za \in Z$ with $\za \neq \zamax$, we have $g(\zamax) \equiv
\gmax > g(\za)$, and ii) the determinant of the $(d-1)\times(d-1)$
matrix of second derivatives of $g$ is nonzero at $\zamax$, then
there exists a constant $D$ such that for all sufficiently large
$n$
\[
    S \,<\, D \times \gmax^n \enspace .
\]
\end{lemma}

Applying Lemma~\ref{lem:cris} to the sum in
Lemma~\ref{lem:prob_calc} we see that we need to maximize
\begin{equation}\label{eq:gr}
g_r(\alpha,\beta,\gamma)
 = \frac{\displaystyle{\left(
    \frac{p(\alpha,\beta,\gamma)}{q^2} \right)^r}}
    {4  \,\gamma^\gamma
        \,(\alpha-\gamma)^{\alpha-\gamma}
        \,(\beta-\gamma)^{\beta-\gamma}
        \,(1-\alpha-\beta+\gamma)^{1-\alpha-\beta+\gamma}}
          \enspace ,
\end{equation}
where for convenience we defined $g_r$ in terms of
$\alpha,\beta,\gamma$ instead of $\zeta_1,\zeta_2,\zeta_3$. We
will show that $g_r$ has a unique maximum at
\[ \zomax =(1/2,1/2,1/4)  \enspace .\]
\begin{lemma}
\label{lem:hyp_gleas} For every $\eps > 0$, and all $k \geq
k_0(\eps)$ if
\[
r \leq   2^{k-1} \ln 2 - \frac{\ln 2}{2} - \frac{1+\eps}{2}
\]
then $ g_r(\zomax) = 1 > g_r(\za)$ for all $\za \in Z$ with $\za
\neq \zomax$. Moreover, the determinant of the matrix of second
derivatives of $g_r$ at $\zomax$ is nonzero.
\end{lemma}
Therefore, for all $r,k,\eps$ as in Lemma~\ref{lem:hyp_gleas}
\[
    \frac{\ex[X^2]}{\ex[X]^2}  \; < \; \frac{D}{A^2} \times
    g_r(\zomax)^n \;=\; D/A^2 \enspace ,
\]
completing the proof of Theorem~\ref{thm:hyp} modulo
Lemmata~\ref{lem:prob_calc}, \ref{lem:cris} and
\ref{lem:hyp_gleas}.\medskip

The proof of Lemma~\ref{lem:prob_calc} is a straightforward
probabilistic calculation. The proof of Lemma~\ref{lem:cris} is
somewhat technical but follows standard asymptotic methods. To
prove Lemma~\ref{lem:hyp_gleas} we will rely very heavily on
Lemma~\ref{gleas}.  In particular, we will show that all local
maxima of $g_r$ occur within a one-dimensional subspace, in which
$g_r$ coincides with the function $\gnae$ of~\eqref{eq:common}.
Specifically, we prove
\begin{lemma} \label{lem:abhalf}
If $(\alpha,\beta,\gamma)$ is a local extremum of $g_r$, then $\alpha
= \beta = 1/2$.
\end{lemma}
This reduces our problem to the one-dimensional maximization for
NAE $k$-SAT, allowing us to easily prove
Lemma~\ref{lem:hyp_gleas}.\smallskip

\noindent{\bf Proof of Lemma~\ref{lem:hyp_gleas}.} Observe that
\begin{eqnarray*}
    g_r(1/2,1/2,\gamma) =  \gnae(2\gamma) \enspace ,
\end{eqnarray*}
where $\gnae$ is the function defined in~\eqref{eq:common} for NAE
$k$-SAT. Thus, the inequality $g_r(\zomax) > g_r(\za)$ for $\za \neq
\zomax$ follows readily from Lemma~\ref{gleas}, giving the first part
of the lemma.

To prove the condition on the determinant of the $3 \times 3$ matrix
of second derivatives, a little arithmetic shows that at $\zomax$ it
is equal to
\[ \frac{256 \,(2^k-2-2kr+2 k^2 r)^2}
                         {2^{4k} q^4}
   \,(4k(k-1) \,r - 2^{2k} \,q^2) \enspace .
\]
Thus, the determinant is negative whenever
\[ 4k(k-1) \,r < 2^{2k} \,q^2 \enspace .\]
For $k=3,4$ this is true for $r<3/2$ and $r<49/12$ respectively,
while for $k \ge 5$ it is true for all $r< \ln 2 \times 2^{k-1}$.
\qed

\section{Proof of Lemma~\ref{lem:prob_calc}}

Recall that $X$ denotes the number of 2-colorings of
$H_k(n,m=rn)$.

\subsection{First moment}
Recall that
\[
q = 1-2^{1-k}\enspace .
\]
The probability that a 2-coloring with $z = \alpha n$ black
vertices and $n-z = (1-\alpha) n$ white vertices makes a random
hyperedge of size $k$ bichromatic is
\[ s(\alpha) = 1 - \alpha^k - (1-\alpha)^k \le q \enspace . \]
Summing over the $2^n$ colorings gives
\[ \ex[X] \,=\, \sum_{z=0}^n {n \choose z} \,s(z/n)^{rn} \enspace .
\]
To bound this sum from below we apply the lower bound of
Lemma~\ref{lem:basic} with $f(\alpha) = s(\alpha)^r$. In
particular, it is easy to see that for all $r>0$
\[
g(\alpha) =
\frac{s(\alpha)^r}{\alpha^{\alpha}(1-\alpha)^{1-\alpha}} =
\frac{(1 - \alpha^k -
(1-\alpha)^k)^r}{\alpha^{\alpha}(1-\alpha)^{1-\alpha}}
\]
is maximized at $\alpha= 1/2$ and that $g(1/2) = 2 q^r$.
Moreover, for any $k > 1$
\[ g''(1/2) = - 8 \,(1-2^{1-k})^{r-1}
 \,\bigl(1 + 2^{1-k} (k (k-1) r - 1) \bigr)
 \,<\, 0 \enspace .
\]
Therefore, we see that there exists a constant $A$ such that
\begin{equation}
\ex[X] \,\ge\, A \times \left( 2 q^r \right)^n \enspace .
\label{eq:firstm}
\end{equation}

\subsection{Second moment}

We first observe that $\ex[X^2]$ equals the expected number of
ordered pairs $S,T$ of 2-partitions of the vertices such that both
$S$ and $T$ are 2-colorings. Suppose that $S$ and $T$ have $\alpha
n$ and $\beta n$ black vertices respectively, while $\gamma n$
vertices are black in both. By inclusion-exclusion a random
hyperedge of size $k$ is bichromatic under both $S$ and $T$ with
probability $p(\alpha, \beta, \gamma)$, \ie
\[1 - \alpha^k - (1-\alpha)^k - \beta^k - (1-\beta)^k
  + \gamma^k + (\alpha-\gamma)^k
  + (\beta-\gamma)^k + (1-\alpha-\beta+\gamma)^k \enspace .
\]
The negative terms above represent the probability that the
hyperedge is mono\-chromatic under either $S$ or $T$, while the
positive terms represent the probability that it is
mono\-chromatic under both (potentially with different colors).
Since the $m=rn$ hyperedges are chosen independently and with
replacement, the probability that all $m=rn$ hyperedges are
bichromatic is $p(\alpha,\beta,\gamma)^{rn}$.

If $z_1$, $z_2$, $z_3$ and $z_4$ vertices are respectively black
in both assignments, black in $S$ and white in $T$, white in $S$
and black in $T$, and white in both, then $\alpha = (z_1+ z_2)/n$,
$\beta = (z_1+ z_3)/n$ and $\gamma = z_1/n$. Thus,
\begin{equation}
 \ex[X^2]  =  \sum
   {n \choose z_1, z_2, z_3, z_4}
   \,p\!\left(\frac{z_1+z_2}{n},
            \,\frac{z_1+z_3}{n}
            \,\frac{z_3}{n} \right)^{rn} \enspace .
\label{eq:secondm}
\end{equation}\qed

\section{Proof of Lemma~\ref{lem:abhalf}}\label{sec:reduce}

We wish to show that at any extremum of $g_r$ we have $\alpha =
\beta = 1/2$.  We start by proving that at any such extremum
$\alpha = \beta$. Note that since, by symmetry, we are free to
flip either or both colorings, we can restrict ourselves to the
case where $\alpha \le 1/2$ and $\gamma \le \alpha/2$.

Let $h(x_1,x_2,x_3,x_4) = \sum_{i} x_i \ln x_i$ denote the entropy
function, and let us define the shorthand $(\p/\p x - \p/\p y) f$
for $\p f/\p x - \p f/\p y$. Also, recall that $q = 1-2^{1-k}$ and
that $p(\alpha, \beta, \gamma) \equiv p$ is
$$
 1 - \alpha^k - (1-\alpha)^k - \beta^k - (1-\beta)^k
  + \gamma^k + (\alpha-\gamma)^k
  + (\beta-\gamma)^k + (1-\alpha-\beta+\gamma)^k .
$$

We will consider the gradient of $\ln g_r$ along a vector that
increases $\alpha$ while decreasing $\beta$. We see
\begin{eqnarray}
\lefteqn{\left( \frac{\p}{\p\alpha} - \frac{\p}{\p\beta} \right)
 \ln g_r(\alpha,\beta,\gamma)} \nonumber \\
& = & \left( \frac{\p}{\p\alpha} - \frac{\p}{\p\beta} \right)
  \Bigl(
   h\bigl(\gamma, \alpha-\gamma, \beta-\gamma, 1-\alpha-\beta+\gamma\bigr)
    \,-\, \ln 4 \,+\, r \,(\ln p - 2 \ln q) \Bigr) \nonumber \\
& = & - \ln (\alpha-\gamma) + \ln (\beta-\gamma)\nonumber \\
& &  + \frac{kr}{p} \,\Bigl( - \alpha^{k-1} + (1-\alpha)^{k-1} +
(\alpha-\gamma)^{k-1}
     + \beta^{k-1} - (1-\beta)^{k-1} - (\beta-\gamma)^{k-1}
   \Bigr)
\nonumber \\
& \equiv & \phi(\alpha) - \phi(\beta) \label{eq:phiab} \enspace ,
\end{eqnarray}
where
\[ \phi(x) \,=\, - \ln (x-\gamma) + \frac{kr}{p}
   \,\left( -x^{k-1} + (1-x)^{k-1} + (x-\gamma)^{k-1} \right)
   \enspace .
\]
Here we regard $p$ as a constant in the definition of $\phi(x)$.

Observe now that if $(\alpha,\beta,\gamma)$ is an extremum of
$g_r$ then it is also an extremum of $\ln g_r$. Therefore, it must
be that $(\p/\p\alpha) \ln g_r = (\p/\p\beta) \ln g_r = 0$ at
$(\alpha,\beta,\gamma)$ which, by~\eqref{eq:phiab}, implies
$\phi(\alpha) = \phi(\beta)$. This, in turn, implies  $\alpha =
\beta$ since $\phi(x)$ is monotonically decreasing in the interval
$\gamma < x < 1$:
\[ \frac{{\rm d}\phi}{\dx} \,=\, - \frac{1}{x-\gamma}
   - \frac{k(k-1)r}{p}
   \,\Bigl( x^{k-2} + (1-x)^{k-2} - (x-\gamma)^{k-2} \Bigr)
   \,<\, 0 \enspace .
\]

Next we wish to show that in fact $\alpha = \beta = 1/2$.  Setting
$\alpha = \beta$, we consider the gradient of $\ln g_r$ along a
vector that increases $\alpha$ and $\gamma$ simultaneously (using
a similar shorthand for $\p g/\p\alpha + \p g/\p\gamma$):
\begin{eqnarray*}
\lefteqn{\left(\frac{\p}{\p\alpha} + \frac{\p}{\p\gamma} \right) \ln g_r(\alpha,\alpha,\gamma)} \\
& = & \left( \frac{\p}{\p\alpha} + \frac{\p}{\p\gamma} \right)
    \Bigl( h\bigl(\gamma, \alpha-\gamma, \alpha-\gamma, 1-2 \alpha+\gamma\bigr) \,-\, \ln 4 \,+\, r \,(\ln p - 2 \ln q) \Bigr) \nonumber \\
& = & -\ln \gamma + \ln (1-2\alpha+\gamma) \\
& & + \frac{kr}{p}
    \left(- 2 \alpha^{k-1} + 2 (1-\alpha)^{k-1} - (1-2\alpha+\gamma)^{k-1} + \gamma^{k-1} \right) \\
& \equiv & \psi(\alpha) \enspace .
\end{eqnarray*}
Clearly, $\psi(\alpha) = 0$ when $\alpha = 1/2$.  To show that
$1/2$ is the only such $\alpha$, we show that $\psi$ decreases
monotonically with $\alpha$ by showing that if $0< \alpha < 1/2$
and $\gamma \leq \alpha/2$, all three terms below are negative for
$k \geq 3$.
\begin{eqnarray*}
\frac{\p\psi}{\p\alpha} & = &
\, - \frac{2}{1-2\alpha+\gamma}\\
    & & +\,  \frac{2k(k-1)r}{p} \left( -\alpha^{k-2} - (1-\alpha)^{k-2} + (1-2\alpha+\gamma)^{k-2} \right) \\
    & & +\,  \frac{2 k^2 r}{p^2} \times \left( \gamma^{k-1} - (1-2\alpha+\gamma)^{k-1} - 2\alpha^{k-1} + 2(1-\alpha)^{k-1} \right) \\
    & & \;\;\;\; \times \, \left( - (\alpha-\gamma)^{k-1} + \alpha^{k-1}
          + (1-2\alpha+\gamma)^{k-1} - (1-\alpha)^{k-1} \right) \enspace.
\end{eqnarray*}

The first and second terms are negative since $1-\alpha >
1-2\alpha+\gamma > 0$, implying $(1-\alpha)^{k-2} >
(1-2\alpha+\gamma)^{k-2}$.  The second factor of the third term is
positive since $f(z)=z^{k-1}$ is convex and $(1-\alpha)-\alpha =
(1-2\alpha+\gamma)-\gamma$ (the factor of 2 on the last two terms
only helps us since $1-\alpha \ge \alpha$).  Similarly, the third
factor is negative since $(1-\alpha)-(1-2\alpha+\gamma) = \alpha -
\gamma \ge \alpha-(\alpha-\gamma) = \gamma$.

Thus, $\p\psi/\p\alpha < 0$ and $\alpha = 1/2$ is the unique
solution to $\psi(\alpha) = 0$. Therefore, if
$(\alpha,\alpha,\gamma)$ is an extremum of $g_r$ we must have
$\alpha = 1/2$.  \qed

\section{Proof of Lemma~\ref{lem:cris}}\label{sec:cris}

\subsection{Preliminaries}
We will use the following form of Stirling's approximation for
$n!$, valid for $n > 0$
\begin{equation}
   \sqrt{2 \pi n} \,n^n \,\e^{-n} \,\left(1+\frac{1}{12 n} \right)
   \,<\, n! \,<\,
   \sqrt{2 \pi n} \,n^n \,\e^{-n} \,\left(1+\frac{1}{6 n} \right)
   \enspace .
\label{eq:stirling}
\end{equation}
We will also use the following crude lower bound for $n!$, valid
for $n \ge 0$
\begin{equation}\label{eq:crude_stir}
n! \ge (n/ \e)^{n} \enspace ,
\end{equation}
using the convention $0^0 \equiv 1$.\medskip

\noindent Let $z_1,\ldots,z_{d}$ be such that $\sum_{i=1}^{d} z_i
= n$. Let $\zeta_i = z_i/n$. Let $\za =
(\zeta_1,\ldots,\zeta_{d-1})$.
\begin{itemize}
\item
If $z_i > 0$ for all $i$, then using the upper and lower bounds
of~\eqref{eq:stirling} for $n!$ and $z_i!$ respectively, and
reducing the denominator further by changing the factor $1+1/(12
z_i)$ to $1$ for $i \ne 1$, we get
\begin{eqnarray}
{n \choose z_1, \cdots, z_d} & < & (2 \pi n)^{-(d-1)/2}
     \,\left( \prod_{i=1}^d \zeta_i^{\,-1/2} \right)
     \,\left( \prod_{i=1}^d \zeta_i^{\,-\zeta_i} \right)^n
     \times \frac{1 + 1/(6n)}{1 + 1/(12z_1)}  \nonumber \\
& \le & (2 \pi n)^{-(d-1)/2}
     \,\left( \prod_{i=1}^d \zeta_i^{\,-1/2} \right)
     \,\left( \prod_{i=1}^d \zeta_i^{\,-\zeta_i} \right)^n
     \enspace , \label{eq:pirogi}
\end{eqnarray}
where for~\eqref{eq:pirogi} we assumed w.l.o.g.\ that $z_1 \le
n/2$. Thus,
\begin{equation}
   {n \choose z_1, \cdots, z_d} \,f(z_1/n,\ldots,z_{d-1}/n)^n \,\le\,
   (2 \pi n)^{-(d-1)/2}
   \,\left( \prod_{i=1}^d \zeta_i^{\,-1/2} \right)
   \,g(\za)^n \enspace .
\label{eq:careful}
\end{equation}
\item
For any $z_i \geq 0$, the upper bound of~\eqref{eq:stirling}
and~\eqref{eq:crude_stir} give
\[
{n \choose z_1, \cdots, z_d}
   \,<\, \frac{7}{6} \,\sqrt{2 \pi n}
        \,\left( \prod_{i=1}^d \zeta_i^{\,-\zeta_i} \right)^n
\enspace ,
\]
implying a cruder bound
\begin{equation}
   {n \choose z_1, \cdots, z_d} \, f(z_1/n,\ldots,z_{d-1}/n)^n \,\le\,
   \frac{7}{6} \,\sqrt{2 \pi n} \,g(\za)^n \enspace .
\label{eq:crude}
\end{equation}
\end{itemize}

\subsection{The main proof}

Our approach is a crude form of the Laplace method for asymptotic
integrals~\cite{debruijn} which amounts to approximating functions
near their peak as Gaussians.

We wish to approximate $g(\za)$ in the vicinity of $\zamax$.  We
will do this by Taylor expanding $\ln g$, which is analytic since
$g$ is analytic and positive.  Since $\ln g$ increases
monotonically with $g$, both $g$ and $\ln g$ are maximized at
$\zamax$.  Furthermore, at $\zamax$ the matrix of second
derivatives of $\ln g$ is that of $g$ divided by a constant, since
\[ \left. \frac{\p^2 \ln g}{\p \zeta_i \,\p
\zeta_j}\right|_{\za=\zamax}\!\!\!
   = \, \frac{1}{\gmax} \frac{\p^2 g}{\p \zeta_i \,\p \zeta_j}
   - \frac{1}{\gmax^2} \frac{\p g}{\p \zeta_i} \frac{\p g}{\p \zeta_j}
\]
and at $\zamax$ the first derivatives of $g$ are all zero.
Therefore, if the matrix of second derivatives of $g$ at $\zamax$
has nonzero determinant, so does the matrix of the second
derivatives of $\ln g$.

Note now that since the matrix of second derivatives is by
definition symmetric, it can be diagonalized, and its determinant
is the product of its eigenvalues.  Therefore, if its determinant
is nonzero, all its eigenvalues are smaller than some
$\lambda_{\max} < 0$. Thus, Taylor expansion around $\zamax$ gives
\[ \ln g(\za) \,\le\, \ln \gmax
   \,+\, \frac{1}{2}\lambda_{\max} \,|\za-\zamax|^2
   \,+\, O(|\za-\zamax|^3)
\]
or, exponentiating to obtain $g$,
\[ g(\za) \,\le\, \gmax
   \,\exp\left( \frac{1}{2}\lambda_{\max} |\za-\zamax|^2 \right)
   \times \bigl(1 + O(|\za-\zamax|^3) \bigr) \enspace .
\]
Therefore, there is a ball of radius $\rho > 0$ around $\zamax$
and constants $Y > 0$ and $g_* < \gmax$ such that
\begin{eqnarray}
\mbox{If } |\za-\zamax| \le \rho,\;\;
  g(\za) & \le & \gmax \,\exp\left( -Y \,|\za-\zamax|^2 \right)
  \label{eq:insidebd} \enspace , \\
  \mbox{If } |\za-\zamax| > \rho,\;\; g(\za) & \le & g_* \enspace
.
  \label{eq:outsidebd}
\end{eqnarray}

We will separate $S$ into two sums, one inside the ball and one
outside:
\[\sum_{\za \in Z:\, |\za-\zamax| \le \rho}
       {n \choose \zeta_1 n, \cdots, \zeta_d n} \,f(\za)^n
 \;\;\;+\; \sum_{\za \in Z:\, |\za-\zamax| > \rho}
       {n \choose \zeta_1 n, \cdots, \zeta_d n} \,f(\za)^n
       \enspace .
\]\smallskip

For the terms inside the ball, first note that if $|\za-\zamax|
\le \rho$ then
\[ \prod_{i=1}^d \zeta_i^{\,-1/2} \le W \mbox{ where }
   W = \left( \min_i \zeta_{\max,i} - \rho \right)^{-d/2} \enspace .
\]
Then, since $|\zeta-\zamax|^2 = \sum_{i=1}^{d-1} (\zeta_i -
\zeta_{\max,i})^2$, using~\eqref{eq:careful} and
\eqref{eq:insidebd} we have
\begin{eqnarray*}
\lefteqn{\sum_{\za \in Z:\, |\za-\zamax| \le \rho}
       {n \choose \zeta_1 n, \cdots, \zeta_d n} \,f(\za)^n} & & \\
   & \le & (2 \pi n)^{-(d-1)/2} \,W \,\gmax^n  \quad \times
       \sum_{z_1,\cdots,z_{d-1} = -\infty}^\infty
       \exp\left( -nY \,\sum_{i=1}^{d-1} (\zeta_i - \zeta_{\max,i})^2
          \right) \\
   & = & (2 \pi n)^{-(d-1)/2} \,W \,\gmax^n  \times \;
         \prod_{i=1}^{d-1} \left( \sum_{z_i=-\infty}^\infty
         \exp\left( -nY \,(z_i/n - \zeta_{\max,i})^2 \right) \right) \enspace .
\end{eqnarray*}
Now if a function $\phi(z)$ has a single peak, on either side of
which it is monotonic, we can replace its sum with its integral
with an additive error at most twice its largest term:
\[ \left| \sum_{z=-\infty}^\infty \phi(z)
   - \int_{-\infty}^\infty \phi(z) \,\dz \,\right| \le 2 \max_z \phi(z)
\]
and so
\begin{eqnarray*}
\sum_{z_i=-\infty}^\infty \exp\left( -nY (z_i/n -
\zeta_{\max,i})^2 \right) & \le & 2 + \int_{-\infty}^\infty
   \exp\left(-nY (z_i/n - \zeta_{\max,i})^2 \right) \,\dz \\
& = & \sqrt{\pi n / Y} + 2 < \sqrt{2 \pi n / Y}
\end{eqnarray*}
where the last inequality holds for sufficiently large $n$.
Multiplying these $d-1$ sums together gives
\begin{equation}
\sum_{\za \in Z:\, |\za-\zamax| \le \rho}
       {n \choose \zeta_1 n, \cdots, \zeta_d n} \,f(\za)^n
   \,\le\, W Y^{-(d-1)/2}\,\gmax^n \enspace .
\label{eq:insidefinal}
\end{equation}

Outside the ball, we use~\eqref{eq:crude}, \eqref{eq:outsidebd}
and the fact that the entire sum has at most $n^{d-1}$ terms to
write
\begin{equation}
   \sum_{\za \in Z:\, |\za-\zamax| > \rho}
        {n \choose \zeta_1 n, \cdots, \zeta_d n} \,f(\za)^n
   \,\le\, n^{d-1} \times \frac{7}{6} \,\sqrt{2 \pi n} \,g_*^n
   \,<\, \gmax^n
\label{eq:outsidefinal}
\end{equation}
where the last inequality holds for sufficiently large $n$.
Combining~\eqref{eq:outsidefinal} and \eqref{eq:insidefinal} gives
\[ S < (W Y^{-(d-1)/2} + 1) \,\gmax^n \equiv D \times \gmax^n \]
which completes the proof.  (We note that the constant $D$ can be
optimized by replacing our sums by integrals and using Laplace's
method \cite{ksat,debruijn}.)

\section{Conclusions}

We have shown that the second moment method yields a very sharp
estimate of the threshold for hypergraph 2-colorability.  It allows us
not only to close the asymptotic gap between the previously known
bounds but, in fact, to get the threshold within a small additive
constant. Yet:\medskip

\noindent $\bullet$ While the second moment method tells us that
\whp\ an exponential number of 2-colorings exist for $r=
\mathrm{\Theta}(2^k)$, it tells us nothing about how to find a
single one of them efficiently.  The possibility that such
colorings actually cannot be found efficiently is extremely
intriguing.

\noindent $\bullet$ While we have shown that the second moment
method works really well, we'd be hard pressed to say why. In
particular, we do not have a criterion for determining a
constraint satisfaction problem's amenability to the method. The
fact that the method fails spectacularly for random $k$-SAT
suggests that, perhaps, rather subtle forces are at
play.
\smallskip

Naturally, one can always view the success of the second moment
method in a particular problem as an aposteriori indication that
the satisfying solutions of the problem are ``largely
uncorrelated". This viewpoint, though, is hardly predictive. (Yet,
it might prove useful to the algorithmic question above).

The solution-symmetry shared by NAE $k$-SAT and hypergraph
2-colorability but not by $k$-SAT, \ie the property that the
complement of a solution is also a solution, explains why the
method gives a nonzero lower bound for these two problems (and why
it fails for $k$-SAT). Yet symmetry alone does not explain why the
bound becomes essentially tight as $k$ grows. In any case, we hope
(and, worse, consider it natural) that an appropriate notion of
symmetry is present in many more problems.

\end{document}